\documentclass{article}

\usepackage{hyperref}       % hyperlinks
\usepackage[numbers]{natbib}
\usepackage[final]{nips_2016}

\usepackage[utf8]{inputenc} % allow utf-8 input
\usepackage[T1]{fontenc}    % use 8-bit T1 fonts
\usepackage{url}            % simple URL typesetting
\usepackage{booktabs}       % professional-quality tables
\usepackage{amsfonts}       % blackboard math symbols
\usepackage{amsmath}
\usepackage{bm}
\usepackage{nicefrac}       % compact symbols for 1/2, etc.
\usepackage{microtype}      % microtypography
\usepackage{algorithm}
\usepackage{algorithmic}
\usepackage{graphicx}
\usepackage{color}

\DeclareMathOperator*{\argmax}{arg\,max}

\newcommand{\ff}{\mathbf{f}}
\newcommand{\xx}{\mathbf{x}}
\newcommand{\xopt}{\ensuremath{\xx_\text{opt}}}
\newcommand{\RR}{\mathbb{R}}
\newcommand{\secref}[1]{\hyperref[#1]{Section~\ref{#1}}}

\title{Interactive Preference Learning of  \\ Utility Functions for Multi-Objective Optimization}

\author{
  Ian Dewancker \hspace{5mm}  Michael McCourt \hspace{5mm} Samuel Ainsworth \vspace{3mm} \\
    SigOpt \\
    San Francisco, CA 94108 \\
    \texttt{\{ian, mike, sam\}@sigopt.com} 
}

\begin{document}

\maketitle

\begin{abstract}
  Real-world engineering systems are typically compared and contrasted using multiple metrics.  For practical machine learning systems, performance tuning is often more nuanced than minimizing a single expected loss objective, and it may be more realistically discussed as a multi-objective optimization problem.  
  We propose a novel generative model for scalar-valued utility functions to capture human preferences in a multi-objective optimization setting.   
  We also outline an interactive active learning system that sequentially refines the understanding of stakeholders ideal utility functions using binary preference queries.
  
\end{abstract}

\section{Introduction}
As machine learning systems become more prominent across industries and organizations, it is important that they be tuned so as
to perform as optimally as possible.
One method that has been increasingly popular in identifying the configuration of an optimal machine learning system is Bayesian black-box optimization of the hyperparameter configurations of machine learning models \cite{snoek2012practical} \cite{ThorntonHutter13} \cite{feurer-nips2015}. 
However, most of these techniques require that the objective be a scalar valued function depending on the hyperparamter configuration
$\xx\in\mathcal{X}$;
in the context of machine learning, this objective is often a scalar-valued cross-validated metric $f(\xx)$.
The space of hyperparameter configurations $\mathcal{X}$ is left intentionally ambiguous.

In practice however, the performance of real systems is often more naturally discussed using a vector of competing metrics
\[
	\ff:\mathcal{X}\to\Omega,\qquad \ff(\xx) = \begin{pmatrix}f_1(\xx)\\\vdots\\f_N(\xx)\end{pmatrix},
\]
where $\Omega$ is the space of possible metric values (in this article, we assume $\Omega\in[0, 1]^N$).
This allows for the perspectives of various stakeholders regarding the aspects of an optimal system to be captured.

For general machine learning systems, this competition may involve trade-offs between predictive accuracy and computational cost of a model \cite{hernandez2015predictive} \cite{shah2016pareto}.
The trade-off between precision and recall can be phrased in this format, and
indeed, all classification problems with unbalanced class representation, such as fraud or spam tasks, might be best addressed using a different metric associated with the loss estimate for each class \cite{fieldsend2009optimizing}.

Indeed, it is standard that engineering systems are designed based on varied input from the involved stakeholders; this naturally leads to a discussion of several metrics related to the performance of the system.  We might pose such a problem,
involving the accumulation of competing metrics and finding an optimal configuration, as
\begin{align*}
\xopt =  \argmax_{\xx \in \mathcal{X}} \ u( \ff(\xx)).
\end{align*}
%Here $\argmax$ is ambiguous and often a set of $\xopt$ values define the \emph{Pareto efficient} \cite{Ehrgott05} points
%which may be considered acceptable.
Here, $u:\Omega\to\RR$ denotes a \emph{utility function} which, in a likely implicit fashion, encapsulates a balance in the preferences between stakeholders.

Of course, the implicit nature of this utility function is a significant stumbling block in any optimization process: it may be
simple for developers of a loan prediction model to specify a need to balance the expected false repayment and false default
predictions, but optimally defining their interaction in the business context is more complicated.
This article proposes a model for this utility function in \secref{gen_inst} based on the idea that the function $u$ is formed
by a product of ``individual utilities'' over individual metrics.
We explore the impact of free parameters in the structure of $u$ and then, in \secref{headings} explain
how these free parameters can be appropriately selected through a thought experiment that stakeholders supervising
the machine learning system conduct before attempting to find $\xopt$.
This concept of interactively conducting a multiobjective optimization was discussed in \cite{Belton2008}.
\secref{sec:experimentalresults} presents some manufactured experiments which demonstrate the viability of this
interactive questioning mechanism to accurately reproduce the behavior of a predefined, but unknown, utility function.

%how to efficiently solicit the ideal tradeoff of False positive vs False negative rates from stakeholders involved in the development of a machine learning system 
%Pareto front exploration, drastically cheaper to conduct a initial thought experiment to determine utility preferences. upfront before resources are spent previous work Brochu paper \cite{Abrahamsen97}
%Here we present a method for interactively learning a utility function for a multi-objective optimization problem that captures utility preferences from stakeholders with a simple user interface, posing binary preference queries
%
%Summary work of methods in \cite{Belton2008}

\section{Generative Model for Multi-objective Utility Functions}
\label{gen_inst}

We propose a utility function composed of a product of one dimensional individual
\footnote{The ideal term here would be marginal utility functions (analogous to marginal densities) but the
	term marginal utility has a specific meaning in the context of economics so we prefer the term individual
	utility instead.}
utility functions %  \ | \ \pmb{\alpha}, \pmb{\beta}
\[
	u(\ff(\xx)) = \prod_{i=1}^{N} u_i(f_i(\xx)).
\]
Each of these individual utility functions could take arbitrary structure; we choose to impose the form of
cumulative distribution functions of beta random variables, so that 
\[
	u_i(f_i(\xx);\,\alpha_i,\beta_i) = \int_0^{f_i(\xx)} \frac{t^{\alpha_i-1}(1-t)^{\beta_i-1}}{B(\alpha_i,\beta_i)}\,\textup{d}t,
\]
where $B$ is the beta function and $\alpha_i$ and $\beta_i$ represent two free parameters governing
the shape of the beta density.  We enforce the belief that these parameters are log-normal
with unknown variance,
\begin{align}
	\label{eq:alphabetadistribution}
  \log{(\alpha_i)} \sim \mathcal{N}(\mu_{\alpha_i}, \sigma_{\alpha_i}), \quad
  \log(\beta_i) \sim \mathcal{N}(\mu_{\beta_i}, \sigma_{\beta_i}),
\end{align}
and describe in \secref{sec:marginallikelihood} how these $4N$ $\mu$ and $\sigma$ parameters are chosen.
We also, at times, use a slight abuse of notation,
$u(\ff(\xx)) \equiv u(\ff(\xx); \, \pmb{\alpha}, \pmb{\beta})$, to
suppress the presence of the parameters in the utility.
This structure (the product of distribution functions) has appeared in other literature, although in the context of
adapting non-stationary data for use in stationary Gaussian processes \cite{snoek2014input}.

Some of the existing literature on designing a utility function involves the use of an additive, rather than
multiplicative, combination of marginal utilities \cite{Belton2008}.  We believe
that the multiplicative structure is potentially more suitable for utility functions with a nonlinear structure,
such as those presented in \secref{sec:experimentalresults} involving constraints or the $F$-score style
utilities designed to balance precision and recall.

%To limit the space of utility functions we construct a generative model for monotonic scalar-valued functions.  \cite{chu2005preference}
%Following a similar idea proposed in Swersky \cite{snoek2014input} to better handle non-stationary functions warping  we use a product of cumulative beta distribution to form the basis of our generative model over utility functions.  Why multiplicative instead of additive??
%\begin{align*}
%u(f_1(\xx), \cdots f_N(\xx) \ | \ \pmb{\alpha}, \pmb{\beta}) &= \prod_{i=1}^{N} \int_0^{f_i(x)} \frac{t^{\alpha_i-1}(1-t)^{\beta_i-1}}{B(\alpha_i,\beta_i)}dt  \\
%\log{(\alpha_i)} \sim \mathcal{N}(\mu_{\alpha_i}, \sigma_{\alpha_i}) \ &, \  
%\log(\beta_i) \sim \mathcal{N}(\mu_{\beta_i}, \sigma_{\beta_i}) 
%\end{align*}

For all $\pmb{\alpha}, \pmb{\beta}$, the resulting utility function $u$ will monotonically increase on the domain $[0,1]^N$.
This is a byproduct of the fact that each individual utility is monotonically increasing which enforces the standard assumption
for multicriteria problems that in an ideal setting all metrics should be optimized and thus any increase in one with no
decrease in others is an improvement in utility.
Specifically, for a particular $\pmb{\alpha}$ and $\pmb{\beta}$, the following will hold,
\vspace{1mm}
\begin{align*}
	u(\ff(\xx_A) \ ; \ \pmb{\alpha}, \pmb{\beta}) \geq \, u(\ff(\xx_B) \ ; \ \pmb{\alpha}, \pmb{\beta}),
	\quad \text{if } f_i(\xx_A) \geq\,  f_i(\xx_B) \text{ for all } i \in \{1, ... , N\}.
\end{align*}

\subsection{Marginal Likelihood for Binary Preference Data \label{sec:marginallikelihood}}

The motivation behind the development of this utility model is the need to balance the input from multiple,
possibly competing, stakeholders during the production of a machine learning system with multiple
competing metrics which define success.  As such, it may be difficult or controversial to judge the value
of a specific set of metrics in absolute terms (i.e., for an engineer or product manager to
assign a scalar value associated with the quality of $\ff(\xx)$).

In contrast, it is generally considered a simpler task to compare two sets of metric values and state
which of the two is ``better'' \cite{chu2005preference}.
To reduce cognitive load on our eventual users, we consider exactly this approach, where the only information that stakeholders must provide is a stated preference between proposed vectors of metric values $\ff_A$ and $\ff_B$.

Our mechanism, described in \secref{headings}, solicits this binary preference between the (implicit) utilities of two multi-objective value configurations with the hope of learning an appropriate model $u$ of the form described in \secref{gen_inst}.
Stakeholders may also lack a significant utility preference between two possible configurations, and our model accounts for this perceptual uncertainty by allowing users to report that two configurations are perceived to have equal\footnote{
	In this setting, we allow for users to specify that objective vectors $\ff_1$ and $\ff_2$ have equal utility either because they are perceived to perform equally well or equally poorly.
}
utility.
Thus, we allow for two types of observations from users: pairs of multi-objective values where a clear preference in utility is observed (denoted by $\mathcal{D}_P$) and pairs of configurations where no preference is specified (denoted by $\mathcal{D}_E$),
\vspace{1mm}
\begin{align*}
	\mathcal{D}_P &= \{ (\mathbf{f}^{p_1}_{1}\prec\mathbf{f}^{p_1}_{2}), \ldots, (\mathbf{f}^{p_M}_{1} \prec \mathbf{f}^{p_M}_{2})  \}   \\
	\mathcal{D}_E &= \{ (\mathbf{f}^{e_1}_{1}\prec\succ \mathbf{f}^{e_1}_{2}), \ldots, (\mathbf{f}^{e_L}_{1} \prec \succ \mathbf{f}^{e_L}_{2}) \}
\end{align*}

We quantify this lack of preference by imposing an insensitivity to utilities which differ by too small a margin; this
margin is defined probabilistically by an equivalence distance $u_E  \sim \mathcal{N}(0, \sigma_{E})$, where $\sigma_E$ is
an additional hyperparameter.

We now define a parametrization strategy based on marginal likelihood to help find the $4N+1$ best hyperparameters
\[
	\pmb{\theta} = ( \mu_{\alpha_1}, \sigma_{\alpha_1}, \mu_{\beta_1}, \sigma_{\beta_1}, \mu_{\alpha_2}, \cdots, \sigma_E )
\]
given specific results $\mathcal{D}_P$ and $\mathcal{D}_E$.  We design this parametrization metric to produce
utility functions that better adhere to observed preferences (because they should be seen as more likely).
We define an auxiliary function
\[
	u_{d}(\ff_1,\ff_2 \ ; \ \pmb{\alpha}, \pmb{\beta}  ) =  \ u(\ff_{2} \ ; \ \pmb{\alpha}, \pmb{\beta} ) -  u(\ff_{1} \ ; \ \pmb{\alpha}, \pmb{\beta})
\]
which describes the distance from utilities for given parameters.  The likelihood function, then, is defined
differently for the sets with and without stated preferences:

\begin{align*}
p( \mathcal{D}_P \ | \ \pmb{\theta}) &= \prod_{i=1}^{M} p ( u(\mathbf{f}^{p_i}_{1}) \prec u(\mathbf{f}^{p_i}_{2}) \ | \ \pmb{\theta} )  \\
&= \prod_{i=1}^{M} \iint h( u_{d}(\mathbf{f}^{p_i}_{1},\mathbf{f}^{p_i}_{2} \ ; \  \pmb{\alpha}, \pmb{\beta}  ) ) \ p( \pmb{\alpha}, \pmb{\beta} \ | \ \pmb{\theta} ) \ d \pmb{\beta} \ d \pmb{\alpha}  \\
p( \mathcal{D}_E \ | \ \pmb{\theta}) &= \prod_{j=1}^{L} p ( u(\mathbf{f}^{e_j}_{1}) \prec \succ u(\mathbf{f}^{e_j}_{2}) \ | \ \pmb{\theta} ) \\
&= \prod_{j=1}^{L} \iint 2 \ p( u_E \leq -\left|u_{d}(\mathbf{f}^{e_j}_{1},\mathbf{f}^{e_j}_{2} \ ; \ \pmb{\alpha}, \pmb{\beta}  )\right| \ ) \ p( \pmb{\alpha}, \pmb{\beta} \ | \ \pmb{\theta} ) d \pmb{\beta} \ d \pmb{\alpha}
\end{align*}
where $h$ is the Heaviside function.  The quantity associated with $p( \mathcal{D}_E \ | \ \pmb{\theta})$  mimics the computation of the p-value in a 2-tailed statistical test.
No simple analytic form is known for this likelihood (so far) so we estimate it using Monte Carlo techniques.

\section{Active Preference Learning of Utility Functions }
\label{headings}

Given the likelihood defined in \secref{sec:marginallikelihood}, it is feasible to take a large set of $\mathcal{D}_P$
(and possibly $\mathcal{D}_E$) results and optimally fit the model parameters $\pmb{\theta}$ with, e.g.,
sequential model-based optimization (SMBO) \cite{bergstra2011algorithms,hutter2011sequential,Shahriari2015}.
This strategy optimizes an expensive, and possibly black-box, objective function
by suggesting inputs to be evaluated and incorporating the resulting observed objective values into an appropriate surrogate model of the objective.
The updated surrogate model is then used to optimally make the next suggestion by identifying
points believed to most benefit the optimization search.

However, before this model fitting can take place, it is necessary to generate the $\mathcal{D}_P$
and $\mathcal{D}_E$ data.  This process involves soliciting preferences from stakeholders, which is
considered an expensive process (likely much more expensive than the likelihood evaluation, even with
the Monte Carlo estimation).  As such, we propose the use of SMBO to identify the best questions to ask so as to
minimize the demands on the stakeholders.
The approach has previously proved useful for conducting interactive optimizations with a
\emph{human in the loop} where the implicit objective is related to that user's perception.  
In particular, Brochu \cite{brochu2010bayesian} \cite{eric2008active} demonstrated the viability of the idea for
interactively learning optimal parameter configurations in computer graphics settings.

In this setting, we no longer seek the optimum of some latent objective; indeed for every possible utility function,
that optimal configuration is always known (simply maximizing each individual metric value).
The aim, instead, is to select preference queries that would help to resolve \emph{uncertainty} about the user's latent utility function.
This process is outlined below in Algorithm~\ref{alg:active_algorithm}.  Initial preference data is collected from the user using no information about the utility model.  Each subsequent preference query begins with fitting the utility model to all the preference data observed so far.
An acquisition function $a$ (discussed in \secref{sec:entropysearch}) that ranks pairs of configurations is then maximized to determine which two configurations to present to the user.  The datasets are updated with the preference information and the process continues.       

\begin{algorithm}[H]
	\caption{Active Utility Function Preference Learning\label{alg:active_algorithm}}
	\begin{algorithmic}
		\STATE {\bfseries Input:}  $\Omega$ 
		\STATE $\mathcal{D}_P, \mathcal{D}_E  \gets \textsc{InitUserPrefs}( \Omega )$
		\FOR{$i \gets 1$ {\bfseries to} $T$}
		\STATE $ \pmb{\theta}_{MLE} \gets \argmax_{\pmb{\theta}} \, p(\mathcal{D}_P \ | \ \pmb{\theta}) \ p(\mathcal{D}_E \ | \ \pmb{\theta}) $
		\STATE $ \ff_{A}, \ff_{B} \gets \argmax_{\ff_1,\ff_2 \in \Omega} \, a( \ff_{1},\ff_2 \ ; \ \pmb{\theta}_{MLE} )$  
		\STATE $ p \gets \textsc{GetUserPref}(\ff_{A}, \ff_{B})   \hskip3em \triangleright \text{ Get preference from user }$ (\{ A,B,E \})
		\IF{$p$ == E}
			\STATE $ \mathcal{D}_E \gets \mathcal{D}_E \, \cup \ ( \ff_{A} \prec \succ \ff_{B}  )  $
		\ELSE
			\STATE $ \mathcal{D}_P \gets \mathcal{D}_P \, \cup \ ( \ff_{o}  \prec \ff_{p}   )  $
		\ENDIF
		\ENDFOR
	\end{algorithmic}
\end{algorithm}

\subsection{Single and Pairwise Maximum Entropy Search \label{sec:entropysearch}}

As described in the previous section, we need an acquisition function $a$ with which to propel the SMBO.
The goal of that acquisition function is to design questions that will reduce our uncertainty in the MLE
of the utility function model.
Entropy based search policies have been used in the literature for selecting instances for labeling in active learning settings for machine learning models \cite{settles2010active}.  In our context, the task is to decide which pair of configurations would give us the most information about our utility function model.

We define an entropy-like condition to power the search for the two metric configurations whose difference in utility
has the greatest uncertainty.  Under a Gaussian assumption, the
entropy of a random variable is a monotonic function of its variance \cite{settles2010active}, so we can use the empirical variance of samples of the utility differences to approximate the entropy,
\begin{align*} 
  a( \ff_{1},\ff_2 \ ; \ \pmb{\theta}_\text{MLE} )
					%   &= H( u_{d}(\ff_1,\ff_2 \ | \ \ \pmb{\theta}_{MLE}  )) \\
                       &= \text{Var}(u_{d}(\ff_1,\ff_2 \ | \ \pmb{\alpha}_\text{MLE} ,\ \pmb{\beta}_\text{MLE} )),
%                       \qquad \pmb{\alpha, \beta} \sim p( \pmb{\alpha, \beta} \ | \ \pmb{\theta}_{MLE}).
\end{align*}
where $\pmb{\alpha}_\text{MLE},\pmb{\beta}_\text{MLE}$ follow the distribution defined in
\eqref{eq:alphabetadistribution} with hyperparameters from $\pmb{\theta}_\text{MLE}$.
This formulation, labeled ``pairwise entropy'', assumes we are searching for two new configurations for the user to compare each iteration.
An alternative strategy is to keep the preferred configuration $\ff_p$ and fix $\ff_1$ to take that value in the next search, varying only $\ff_2$.  This search strategy is labeled ``single entropy''.
The distinction is pertinent in \secref{sec:experimentalresults}.

It is important to note that we enforce that the configuration pairs $\ff_1, \ff_2$ presented to users always exhibit the property $\exists i, \exists j \ \ff_1[i] > \ff_2[i], \ff_2[j] > \ff_1[j] $, so that we are not merely asking users to compare configurations with an already assumed preference based on our utility being monotonic.

\section{Interactive Tool for Utility Preference Queries  }
\label{others}

To realize the full benefit of incorporating many perspectives on the optimality of a system, we outline a simple interface for our utility preference solicitation system that is easy to understand even for groups of users having a broad range of technical sophistication.  

\subsection{Multi-objective Utility Comparison Cards}

\begin{figure}[H]
	\centering
	\includegraphics[width=0.85\linewidth]{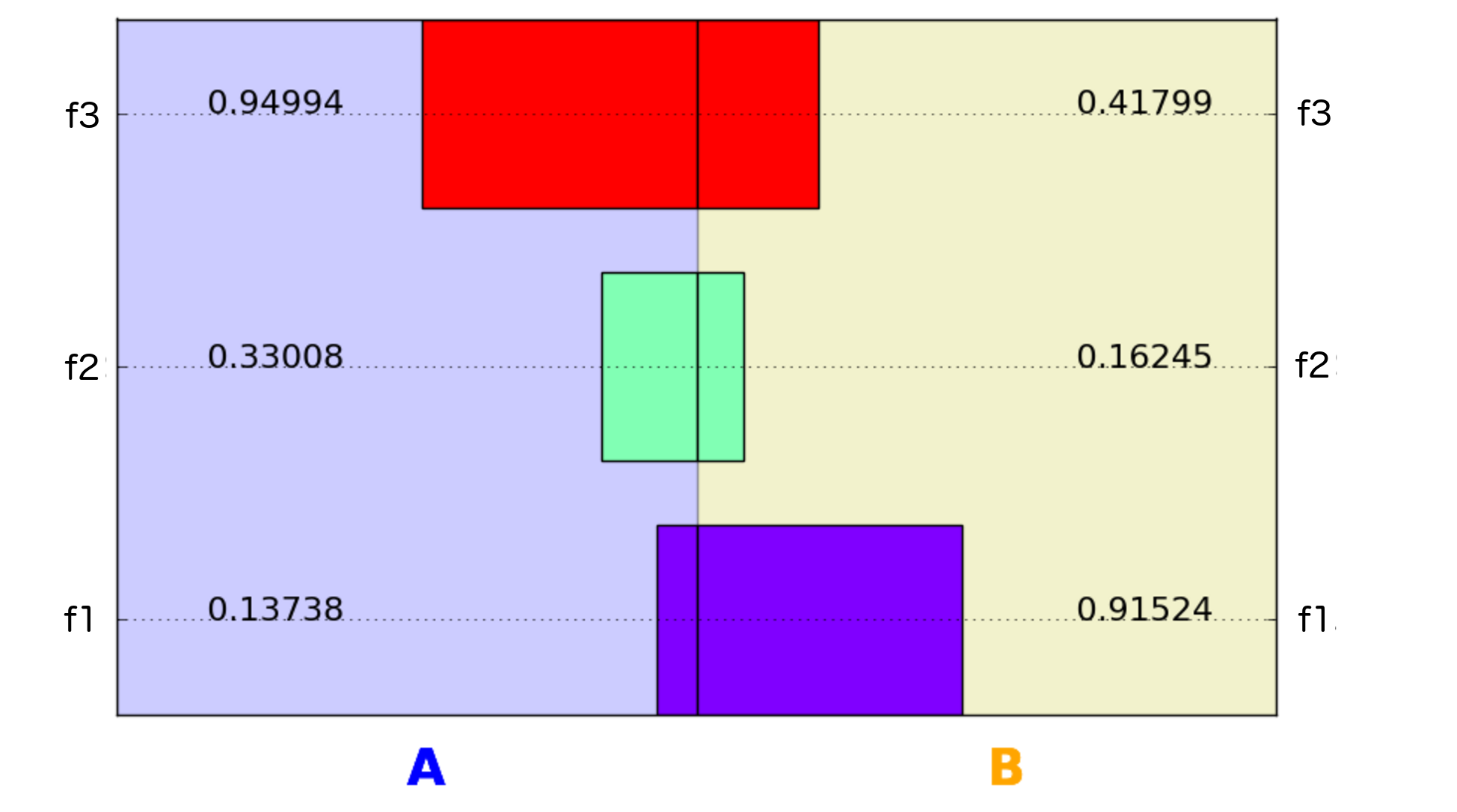}
	\caption{A sample comparison card for preference solicitation showing 3 metrics ($f_1, f_2, f_3$) in two configurations A and B.  Users are asked if they believe utility of configuration A or B to be higher.  Users can also specify that they perceive the utilities of the configurations as equal}
	\label{fig:pref_card}
\end{figure}

To visualize two multi-objective value configurations, our system uses a simple back-to-back bar chart as shown in Figure~\ref{fig:pref_card}.  Users are asked to select which configuration they perceive as having higher utility.   It is worth noting that alternate visualizations could be used in place of the one outlined here, including ones that do not directly expose the numeric value of the underlying metrics but rely on more system-specific visual summaries.

\subsection{Visualization of Utility Functions  }

\begin{figure}[H]
	\centering
	\includegraphics[width=\linewidth]{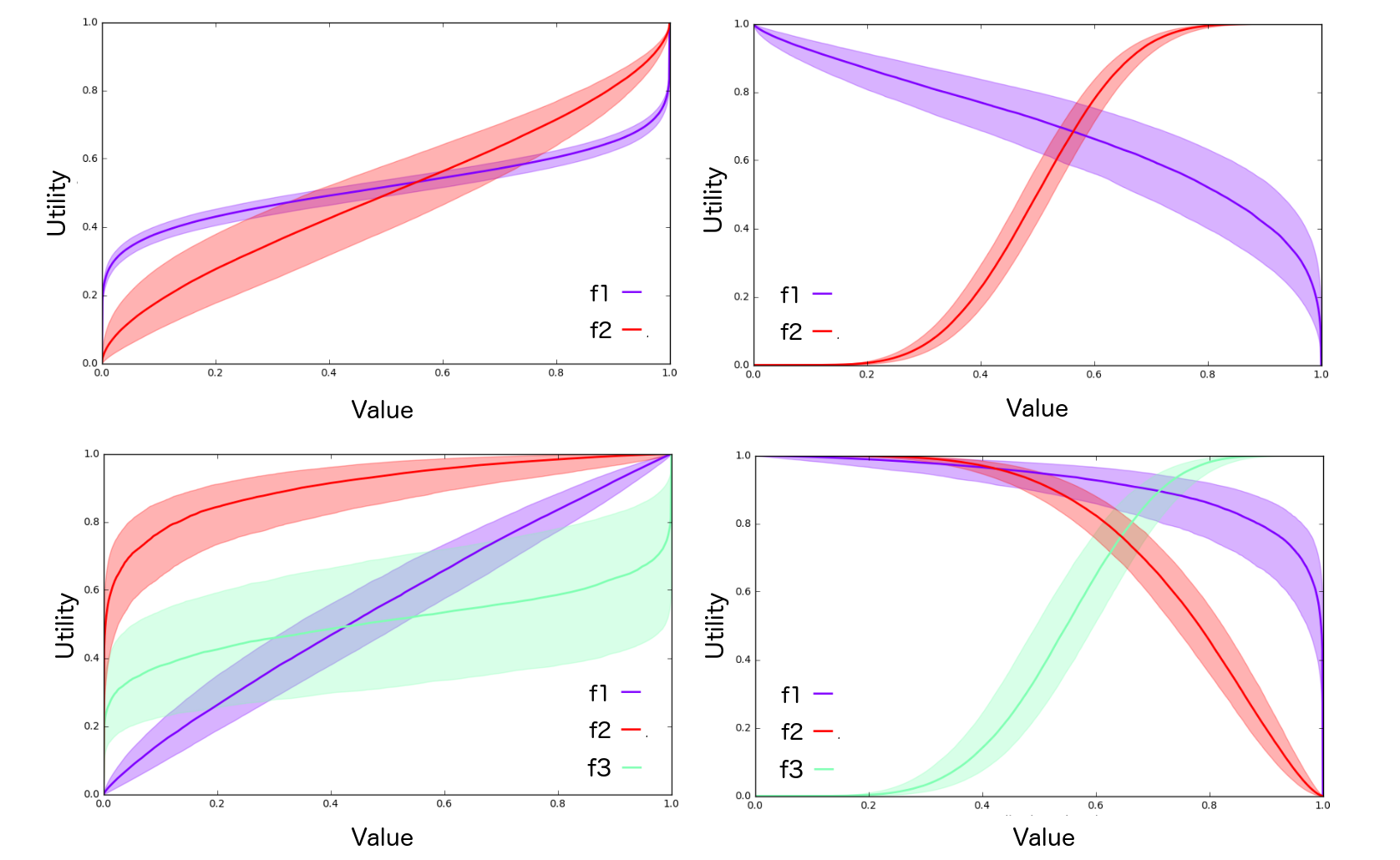}
	\caption{Plots of learned independent utilities with mean and interquartiles.   Starting from top left and proceeding clockise, the test utility functions were : 1. max $ f_1 + 2f_2$  2. min $ f_1 \ \text{s.t.} \ f_2 > 0.6$ \newline 3. max $5f_1 + 2f_2 + f_3$  4. min $f_1 \ \text{s.t.} \ f_2 < 0.2, f_3 > 0.6 $}
	\label{fig:intro_viz}
\end{figure}

\begin{figure}[ht]
	\centering
	\includegraphics[width=\linewidth]{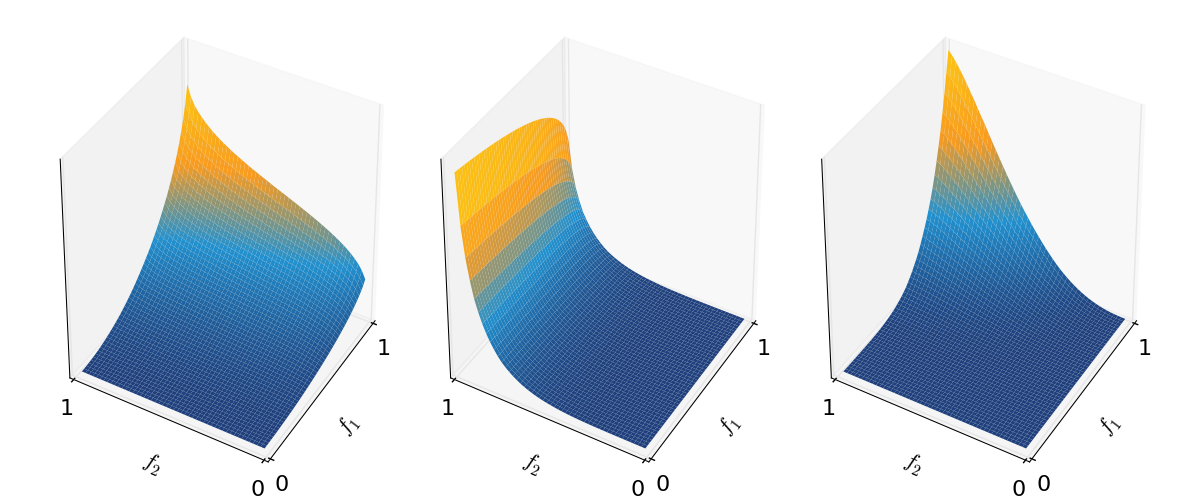}
	\caption{Examples of the learned utility functions.
		\textit{left}: max $f_1 + 2 f_2$,
		\textit{center}: min $f_1 \text{ s.t. } f_2 > 0.6$,
		\textit{right}: max  $5\, f_1 f_2\,/\,(4f_1+f_2) $.
	}
	\label{fig:sampleutilitysurfaces}
\end{figure}

Since the learned utility is a product of $N$ cumulative distribution functions, these distributions can 
be independently visualized to give the user an intuition about the components of the full product utility function.  The MLE of the $4N$ hyperparameters governing the distribution of the $\pmb{\alpha, \beta}$ terms can serve users as a direct introspection mechanism into the uncertainty associated with the learned individual utilities.  Figure~\ref{fig:intro_viz} shows an example of these independent utility plots for several test utility functions with their mean and interquartile range highlighted;
Figure~\ref{fig:sampleutilitysurfaces} shows the actual utility learned from three sample implicit utilities.

Figure \ref{fig:intro_viz} provides some confidence in the learned joint utility functions.  For example, consider the top right plot (2), corresponding to the test utility : min $ f_1 \ \text{s.t.} \ f_2 > 0.6$.  We can see that the model has attempted to learn the threshold constraint of $f_2 > 0.6$.  We see that the individual utility for $f_2$ has a sharp, non-linear spike around 0.6.  In the full utility function product then, configurations with $f_2 < 0.6$ will be zeroed out and those with $f_2 > 0.6$, the utility function will take on the values of $f_1$.  We can also see that the individual utility plot of $f_1$ is a mostly linear looking monotonically decreasing function, maximized when $f_1 = 0$ and minimized when $f_1 = 1$ which corresponds nicely to a utility function for a metric we aim to minimize.  We allow in the specification of the model for metrics to defined as optimally minimized or maximized.  The individual utility for minimization metrics is defined as the survival function $(1-u_i(f_i(\xx))$ of the beta cumulative distribution.       

\section{Experimental Results \label{sec:experimentalresults}}
\label{exp_results}

A series of experiments were conducted to evaluate the performance of the learned utility function model using several active search policy variations.  Explicit test utility functions were used to simulate implicit human utility functions.   We include 
two test utilities that correspond to the $F_1$ and $F_2$ scores commonly used in machine learning and information retrieval \cite{Rijsbergen}.  We also include test utility functions that incorporate threshold constraints.  These threshold constrained tests were simulated in the following way: if both configurations violated the constraints, the configurations were reported with equal perceived utility.

 A hold-out set of 10,000 random multi-objective configurations were generated for each test function and the Kendall rank correlation coefficient was used to quantify the ordinal association between the test utility function values and the learned utility function values for all 10,000 configurations.  Since the region of the feasible solutions is not known \emph{a priori}, our learned utility must strive to ensure that the utility is consistent over the entire domain.   

Each method was allowed $10N$ binary preference queries  
where $N$ is the number of metrics.  The active search strategies were initialized with $5N$ randomly selected configurations.  We report the average Kendall correlation score after 5 runs using each search algorithm on the same generated hold-out set.

\begin{table}[H]
	\caption{Kendall-Tau Correlation using Different Search Policies}
	\label{results-table}
	\centering
	\begin{tabular}{l c c c }
		\cmidrule{1-4} 
		Test Utility Function     & Rnd Search & Single Entropy  & Pair Entropy \\
		\midrule
		max $f_1 + 2 f_2$  & \bf{0.8756}  &  0.8542 & 0.8618   \\
		max $f_1 + 10 f_2$  & 0.9422  & 0.9448 & \bf{0.9615}   \\
		min $f_1 \text{ s.t. } f_2 > 0.6$  & 0.6507 & 0.6805 & \bf{0.6893}      \\
		max $2\, f_1 f_2\,/\,(f_1+f_2) $       & 0.8844 & 0.9028 & \bf{0.9039}  \\
%		max  $2 \frac{f_1 f_2}{f_1+f_2} $       & 0.8844 & 0.9028 & \bf{0.9039}  \\
		max  $5\, f_1 f_2\,/\,(4f_1+f_2) $       & 0.8949 & 0.8950 & \bf{0.9120}  \\
%		max  $5 \frac{f_1 f_2}{4f_1+f_2} $       & 0.8949 & 0.8950 & \bf{0.9120}  \\
		max  $f_1 + 2f_2 + f_3 $ & \bf{0.8490} & 0.8018 & 0.7805  \\
		max  $5f_1 + 2f_2 + f_3$ & \bf{0.8738} & 0.8516 & 0.8311  \\
		min $f_1$ s.t. $f_2 > 0.6, f_3 < 0.2$ &  0.2949 & 0.3154 & \bf{0.3257}  \\
		max  $2\, (f_1 f_2)\,/\,(f_1+f_2) \text{ s.t. } f_3 > 0.95 $       & 0.2309 & 0.2088 & \bf{0.2648}  \\
%		max  ($2 \frac{f_1 f_2}{f_1+f_2} \text{ s.t. } m_3 > 0.95 $       & 0.2309 & 0.2088 & \bf{0.2648}  \\
		\bottomrule
	\end{tabular}
\end{table}

From these results, our proposed utility model appears to perform well under the various acquisition functions.  The performance of the random search policy is particularly noteworthy as it avoids fitting the utility model each iteration and therefore much less expensive computationally.
The mean performance of the pair entropy search seems to edge out the other methods on most of the examples with the interesting exception of the linear test utility functions.

Future work could involve experiments using real human users on a relevant machine learning system building task.  In addition, further investigations into the utility function model and acquisition function could prove valuable in capturing certain utilities.

\small

\small
\bibliographystyle{plain}
\bibliography{citations}

\end{document}